\documentclass[11pt]{amsart}

\usepackage[T1]{fontenc}
\usepackage{lmodern}
\usepackage{amsmath,amssymb}
\usepackage[expansion=false]{microtype}
\usepackage[a4paper,margin=28mm]{geometry}
\usepackage{xcolor}
\usepackage[hidelinks]{hyperref}
\everymath{\displaystyle}

\newtheorem{theorem}{Theorem}[section]
\newtheorem{proposition}[theorem]{Proposition}
\newtheorem{lemma}[theorem]{Lemma}
\newtheorem{corollary}[theorem]{Corollary}
\numberwithin{equation}{section}

\newcommand{\HH}{\mathcal H}
\newcommand{\KK}{\mathcal K}
\newcommand{\BH}{\mathcal B(\mathcal H)}
\newcommand{\Idem}{\mathfrak E}
\newcommand{\CE}{\mathfrak C_{\mathfrak E}}
\newcommand{\DE}{\mathfrak D_{\mathfrak E}}
\newcommand{\ran}{\operatorname{ran}}
\newcommand{\diag}{\operatorname{diag}}
\newcommand{\spancl}{\overline{\operatorname{span}}}
\newcommand{\ueq}{\simeq}

\title[Idempotent commutators and nilpotent similarity]
{O\lowercase{ptimal lengths of commutators of idempotents and similarity of nilpotent operators to their negatives}}
\author{M\lowercase{ohamed} A\lowercase{mine} A\lowercase{ouichaoui}}
\address{Department of Mathematics, Faculty of Sciences of Monastir, University of Monastir, 5019 Monastir, Tunisia}
\email{aouichaoui.mohamedamine@fsm.u-monastir.tn}
\email{amine.aouichaoui@fsm.rnu.tn}
\date{}

\begin{document}

\begin{abstract}
Let $\HH$ be a complex separable infinite-dimensional Hilbert space.
We prove that every operator in $\BH$ is a sum of three commutators
of idempotents and a product of two such commutators. Both numbers
are best possible. The same is true for differences of idempotents,
for operators similar to their negatives, and for operators
anticommuting with an involution. Each of these four classes therefore
has additive length three and multiplicative length two. We also
construct nilpotent operators of order three that are similar to
their negatives but cannot be conjugated to their negatives by any
bounded involution. These results answer Question~7.22(a)--(c),
the question in Remark~7.17, and the question about involutions
in Section~2.9 of Marcoux, Radjavi and Zhang,
\emph{J. Funct. Anal.} \textbf{284} (2023), 109854.
\end{abstract}

\maketitle

\section{Introduction and main results}

Throughout this paper, $\HH$ denotes a complex separable infinite-dimensional
Hilbert space, and $\BH$ denotes the algebra of bounded linear operators on $\HH$. We write $I$ for the identity operator and
$[A,B]:=AB-BA$. Following \cite{MRZ}, let
\[
 \Idem=\{E\in\BH:E^2=E\},\qquad
 \CE=\{[E,F]:E,F\in\Idem\},\qquad
 \DE=\{E-F:E,F\in\Idem\}.
\]
For $\Delta\subseteq\BH$ and a positive integer $n$, let
$\sum_n\Delta$ and $\prod_n\Delta$ denote the sets of sums and
products of $n$ elements of $\Delta$. Define
\[
 \ell_+(\Delta)=\min\{n\geq1:\textstyle\displaystyle\sum_n\Delta=\BH\},
 \qquad \text{and} \qquad
 \ell_\times(\Delta)=\min\{n\geq1:\textstyle\displaystyle\prod_n\Delta=\BH\}.
\]
The minima are over positive integers. A length is $\infty$ if no
such integer exists. The factors in a product need not be invertible.
Each class in Theorem~\ref{thm:main} below contains $0$, so shorter sums can
be extended by zero summands.
An involution is an operator $J$ with $J^2=I$; it need not be
self-adjoint. We write $A\sim B$ for similarity and $A\ueq B$ for
unitary equivalence. The four classes of operators in Theorem~\ref{thm:main} are invariant
under unitary conjugation. We use fixed unitary identifications to
regard operators on finite or countable orthogonal sums of copies
of $\HH$ as operators on $\HH$. Write $\HH^{(\infty)}=\bigoplus_{n=1}^{\infty}\HH.$ For $A\in\BH$, write $A^{(\infty)}=\bigoplus_{n=1}^{\infty}A.$

For unrestricted bounded factors, Brown and Pearcy \cite{BrownPearcy}
proved that an operator in $\BH$ is a commutator if and only if it
is not of the form $\lambda I+K$, where $\lambda\ne0$ and $K$ is
compact. See also Halmos \cite[pp.~128--134]{Halmos} for a discussion
of operator commutators. Wu \cite{WuFactorization} surveys
factorizations into prescribed classes, including involutions and
partial isometries. Here the two entries of each commutator are
required to be idempotents, and we determine the least numbers of
summands and factors needed to obtain every operator in $\BH$.

Marcoux, Radjavi and Zhang proved
\[
 \BH=\sum_5\CE=\prod_3\CE=\sum_4\DE,
 \qquad I\notin\sum_2\CE;
\]
see \cite[Propositions~7.7 and~7.18, Theorem~7.16]{MRZ}.
They also proved that $\prod_2\CE$ is norm dense in $\BH$ and asked
whether it equals $\BH$ \cite[Remark~7.17]{MRZ}.
Question~7.22 asks for the least number of differences of idempotents
needed in a sum, whether every thin operator is a product of two such
differences, and the least numbers of summands and factors for operators
similar to their negatives. Here a thin operator is a scalar multiple
of $I$ plus a compact operator.

\begin{theorem}\label{thm:main}
Let $\Delta\subseteq\BH$ be one of the following four classes:
commutators of idempotents, differences of idempotents, operators
similar to their negatives, or operators anticommuting with an
involution. Then,
\[
 \ell_+(\Delta)=3,\qquad \ell_\times(\Delta)=2.
\]
\end{theorem}

This gives $n=3$ in Question~7.22(a), an affirmative answer to
Question~7.22(b), and $m=2$, $n=3$ in Question~7.22(c).
The equality $\BH=\prod_2\CE$ answers the question in Remark~7.17.

Note that the equality $\ell_+(\DE)=3$ also follows from earlier results.
Pearcy and Topping proved that every operator in $\BH$ is a sum of
five idempotents \cite{PT}; see also the proof of
\cite[Theorem~2]{RS}. Apply their theorem to $T+3I$, add a zero
idempotent, and use $E+F-I=E-(I-F)$ to pair the six terms.
The lower bound follows from \cite[Theorem~1]{RS}, which excludes
$-I$ from sums of four idempotents.

For every $T\in\BH$ and every real $c>\|T\|$,
Theorem~\ref{thm:products} gives $A,B\in\CE$ satisfying
\[
 T=AB,\qquad AA^*=B^*B=cI,\qquad \|A\|=\|B\|=\sqrt c.
\]
Fong and Sourour \cite{FongSourour} proved that, if a unitary is not
a compact perturbation of a scalar operator, then every unitary is a
finite product of operators unitarily equivalent to it. For the
bilateral shift of multiplicity one, three factors suffice by
\cite[Theorem~5.3]{Tricycles}, and two suffice by
\cite[Theorem~1.1]{Bicycles}.
Proposition~\ref{prop:unitary-factors} uses the latter theorem to
write every unitary $T$ as $T=AB$, where $A$ and $B$ are unitary
commutators of idempotents. Both commutators are obtained by unitary
conjugation of a fixed pair of idempotents. For this factorization,
the identities above hold with $c=1=\|T\|$. In Proposition~\ref{prop:fixed}, the right factor is independent of
$T$; both factors anticommute with involutions and are differences
of idempotents.

Drnov\v{s}ek, Radjavi and Rosenthal characterized commutators of
idempotents in \cite[Theorem~1]{DRR}. Drnov\v{s}ek gave an equivalent condition in terms of involutions
in \cite[Theorem~2.1]{Drnovsek}. We use this theorem to show that $\diag(a,a,-a,-a)$ is a commutator
of two idempotents over every unital associative algebra whose scalar
field has characteristic different from two. Countable direct sums and
\cite[Lemma~2.6(ii)]{DM} give the additive result.
For the multiplicative result, we write each strict contraction as
$V^*W$, where $V,W$ are pure isometries with infinite-dimensional
$\ker V^*$ and $\ker W^*$. Proposition~\ref{prop:fourfold} then shows
that the multiples of $V^*$ and $W$ used in the factorization belong
to $\CE$.

An operator $T\in\BH$ anticommutes with an involution if and only if
it is a scalar multiple of a commutator of idempotents. Indeed, if $JT=-TJ$ and $J^2=I$, then
\cite[Corollary~2.4]{Drnovsek}, applied to a sufficiently small
positive multiple of $T$ relative to the eigenspace decomposition
of $J$, gives idempotents $P,Q$ and $\delta>0$ such that
$[P,Q]=\delta T$. Conversely, $2P-I$ is an involution that anticommutes with every
scalar multiple of $[P,Q]$. This equivalence is also recorded in
\cite[Proposition~2.2]{Quadratic}.

In \cite[Section~2.9]{MRZ}, the authors also ask whether a nilpotent
operator similar to its negative must be conjugated to its negative
by an involution. We give a negative answer for nilpotent operators
of order three.

\begin{theorem}\label{thm:nilpotent-main}
There are $M\in\BH$ and an invertible $S\in\BH$ such that
\[
 M^3=0\neq M^2,\qquad S^{-1}MS=-M,
\]
but no involution $J\in\BH$ satisfies $JMJ=-M$.
\end{theorem}
\newpage
Section~\ref{sec:nilpotents} allows the upper-right entry in the block
form of \cite[Section~2.6]{MRZ} to vary. Under the hypotheses of
Theorem~\ref{thm:involution-criterion}, a bounded involution exists
exactly when the ratios of the block coefficients are bounded.
The involution need not preserve the summands.
Theorem~\ref{thm:all-schatten} proves Theorem~\ref{thm:nilpotent-main}
with an operator whose singular values are $p$-summable for every
$p>0$. For each $\eta>0$, it gives such operators $M,S$ with
$\max\{\|S\|,\|S^{-1}\|\}\leq1+\eta$.

The paper is organized as follows. Section~2 determines the additive
lengths. Section~3 proves the product representations, including the
unitary case and the factorization with a fixed right factor.
Section~\ref{sec:nilpotents} constructs the nilpotent examples.

\section{Sums of three commutators}

\begin{proposition}\label{prop:fourfold}
Let $\mathcal A$ be a unital associative algebra over a field of
characteristic different from two. For every $a\in\mathcal A$,
$\diag(a,a,-a,-a)$ is a commutator of two idempotents in $M_4(\mathcal A)$.
\end{proposition}

\begin{proof}
Let $1$ denote the identity of $\mathcal A$, and set
\[
 \widehat a=\begin{pmatrix}a&0\\0&a\end{pmatrix},
 \qquad
 r=\begin{pmatrix}0&a^2+\frac14\,1\\1&0\end{pmatrix}.
\]
Then, $r^2=\widehat a^{\,2}+\frac14 I_2$ and
$r\widehat a=\widehat a r$. In $M_4(\mathcal A)$ put
\[
 t=\begin{pmatrix}\widehat a&0\\0&-\widehat a\end{pmatrix},
 \qquad
 s=\begin{pmatrix}r&0\\0&r\end{pmatrix},
 \qquad
 j=\begin{pmatrix}0&I_2\\I_2&0\end{pmatrix}.
\]
These elements satisfy
\[
 j^2=I_4,\qquad jt=-tj,\qquad js=sj,\qquad st=ts,
 \qquad s^2=t^2+\frac14 I_4.
\]
These are the conditions in \cite[Theorem~2.1(iii)]{Drnovsek}.
Their equivalence with condition~\textup{(i)} proves that $t$ is a
commutator of two idempotents.
\end{proof}

We use the alternating direct sums from the proof of
\cite[Proposition~7.18]{MRZ}, together with
Proposition~\ref{prop:fourfold}, to obtain the following statements
for commutators of idempotents.

\begin{corollary}\label{cor:stable}
For every $A\in\BH$,
\[
 A^{(\infty)}\oplus(-A)^{(\infty)}\in\CE,
 \qquad
 A\oplus0\in\sum_2\CE,
\]
where the zero summand acts on an infinite-dimensional space. There is
also a function $f:\BH\to\BH$ such that $A\oplus f(A)\in\CE$ for every
$A\in\BH$.
\end{corollary}

\begin{proof}
Apply Proposition~\ref{prop:fourfold} to $A^{(\infty)}$ and use
$A^{(\infty)}\oplus A^{(\infty)}\ueq A^{(\infty)}$.
Regrouping the copies of $A$ and $-A$ gives the first assertion.

On $\HH^{(\infty)}$, define
\[
 Y_A=\diag(A,-A,A,-A,\ldots),
 \qquad
 Z_A=\diag(0,A,-A,A,-A,\ldots).
\]
The first assertion and a permutation of the summands give
$Y_A\in\CE$. Since $Z_A\ueq0\oplus Y_A$ and
\[
 0\oplus[E,F]=[0\oplus E,0\oplus F],
\]
we also have $Z_A\in\CE$. Their sum is
\[
 Y_A+Z_A=\diag(A,0,0,\ldots)\ueq A\oplus0.
\]
Finally, fix a unitary $W:\HH\to\HH^{(\infty)}$ and put
\[
 f(A)=W^*\diag(-A,A,-A,A,\ldots)W.
\]
Then, $A\oplus f(A)\ueq Y_A\in\CE$.
\end{proof}

We recall the following result of Davidson and Marcoux
\cite[Lemma~2.6(ii)]{DM}.

\begin{lemma}\label{lem:gluing}
Let $f:\BH\to\BH$ be a function, and let $\mathcal M_1,\mathcal M_2$
be unitarily invariant subsets of $\BH$. Suppose that
\[
 A\oplus f(A)\in\mathcal M_1,
 \qquad
 B\oplus0\in\mathcal M_2
 \qquad(A,B\in\BH),
\]
where the zero summand is infinite-dimensional. Then,
$\BH=\mathcal M_1+\mathcal M_2$.
\end{lemma}

Every commutator of idempotents anticommutes with an involution and
is so similar to its negative; see
\cite[Proposition~2.1]{MRZclosures}.
Every operator anticommuting with an involution is also a difference
of idempotents: relative to the topological direct sum of the two
eigenspaces, its diagonal blocks vanish, and the calculation in the
proof of \cite[Proposition~2.2]{MRZclosures} applies.

\begin{theorem}\label{thm:additive}
One has
\[
 \BH=\sum_3\CE=\sum_3\DE,
 \qquad -3I\notin\sum_2\DE.
\]
Every operator is also a sum of three operators anticommuting with
involutions. The identity is not a sum of two operators similar to
their negatives. Hence $\ell_+(\Delta)=3$ for each of the four
classes in Theorem~\ref{thm:main}.
\end{theorem}

\begin{proof}
The classes $\CE$ and $\sum_2\CE$ are unitarily invariant.
Corollary~\ref{cor:stable} gives the hypotheses of
Lemma~\ref{lem:gluing} with $\mathcal M_1=\CE$ and
$\mathcal M_2=\sum_2\CE$. Hence
\[
 \BH=\CE+\sum_2\CE=\sum_3\CE.
\]
The inclusions stated before the theorem give the remaining three
upper bounds.

The spectral argument in the proof of \cite[Proposition~7.7]{MRZ}
excludes $I=A+B$ whenever $A$ and $B$ are similar to their negatives.
It therefore also excludes two commutators of idempotents or two
operators anticommuting with involutions.

By \cite[Theorem~1]{RS}, $-I$ is not a sum of four idempotents. If
\[
 -3I=(E_1-F_1)+(E_2-F_2),\qquad E_j,F_j\in\Idem,
\]
then
\[
 -I=E_1+(I-F_1)+E_2+(I-F_2)
\]
would be such a sum. Thus, $-3I\notin\sum_2\DE$.
\end{proof}

\section{Products of two commutators}

An isometry $V\in\BH$ is called pure if
$\bigcap_{k\geq0}V^k\HH=\{0\}$. It has infinite multiplicity when
$\ker V^*$ is infinite-dimensional.

\begin{theorem}\label{thm:purefactor}
For every $C\in\BH$ with $\|C\|<1$, there are pure isometries
$V,W\in\BH$ such that $C=V^*W$ and both $\ker V^*$ and
$\ker W^*$ are infinite-dimensional.
\end{theorem}

\begin{proof}
Let $D=(I-C^*C)^{1/2}$. Since
\[
 D\geq\sqrt{1-\|C\|^2}\,I>0,
\]
$D$ has a bounded inverse. On
$\KK=\HH\oplus\HH\oplus\HH$, define
\[
 E,F:\HH\longrightarrow\KK,
 \qquad Ex=(x,0,0),\quad Fx=(Cx,Dx,0).
\]
The equality $C^*C+D^2=I$ shows that $E$ and $F$ are isometries,
with $E^*F=C$. We choose a unitary $U:\KK\to\HH$ such that both
$UE$ and $UF$ are pure.

Fix an orthonormal basis $(e_j)_{j\geq1}$ of $\HH$, and put
\[
 H_n=\spancl\{e_j:j\geq n\},
 \qquad M_n=EH_n+FH_n\quad(n\geq1).
\]
The operator
\[
 L:\HH\oplus\HH\longrightarrow\HH\oplus\HH,
 \qquad L(x,y)=(x+Cy,Dy)
\]
is boundedly invertible, with
\[
 L^{-1}(a,b)=(a-CD^{-1}b,D^{-1}b).
\]
Thus,
\[
 M_n=L(H_n\oplus H_n)\oplus\{0\}.
\]
In particular, the $M_n$ are closed and decreasing, and
\[
 M_1=\HH\oplus\HH\oplus\{0\},
 \qquad \bigcap_{n\geq1}M_n=\{0\}.
\]
The subspace $H_{n+1}\oplus H_{n+1}$ has codimension two in
$H_n\oplus H_n$. Since $L$ is invertible, $M_{n+1}$ has codimension
two in $M_n$. Therefore,
\[
 G_n=M_n\cap M_{n+1}^{\perp}
\]
has dimension two. The orthogonal projections onto $M_n$ decrease
strongly to zero. Consequently,
\[
 M_n=\bigoplus_{k\geq n}G_k,
 \qquad
 \KK=M_1^\perp\oplus\bigoplus_{n\geq1}G_n.
\]

Choose an orthonormal basis $a_n,b_n$ of each $G_n$, and an orthonormal
basis $(z_n)_{n\geq1}$ of
$M_1^\perp=\{0\}\oplus\{0\}\oplus\HH$. Define $U:\KK\to\HH$ on
these basis vectors by
\[
 Uz_n=e_{3n-2},
 \qquad Ua_n=e_{3n-1},
 \qquad Ub_n=e_{3n}.
\]
This defines a unitary. Moreover,
\[
 UM_n=\spancl\{e_{3k-1},e_{3k}:k\geq n\}\subseteq H_{n+1}.
\]
Now put $V=UE$ and $W=UF$. These are isometries on $\HH$, and
\[
 V^*W=E^*U^*UF=E^*F=C.
\]
The inclusions $EH_n,FH_n\subseteq M_n$ give
\[
 VH_n\subseteq H_{n+1},
 \qquad WH_n\subseteq H_{n+1}.
\]
Since $H_1=\HH$, it follows by induction that
\[
 V^k\HH\subseteq H_{k+1},
 \qquad W^k\HH\subseteq H_{k+1}\quad(k\geq0).
\]
The intersection of the $H_n$ is zero, so both isometries are pure.
Finally, $\ran V,\ran W\subseteq UM_1$, and so
\[
 U(M_1^\perp)\subseteq\ker V^*\cap\ker W^*.
\]
The subspace $U(M_1^\perp)$ is infinite-dimensional, so both kernels
are infinite-dimensional.
\end{proof}

\begin{lemma}\label{lem:shifts}
Let $V\in\BH$ be a pure isometry with infinite-dimensional $\ker V^*$.
Then, $\lambda V,\lambda V^*\in\CE$ for every $\lambda\in\mathbb C$.
There is also a self-adjoint unitary $G$ such that
$GVG=-V$ and $GV^*G=-V^*$.
\end{lemma}

\begin{proof}
By the Wold decomposition \cite[Chapter~I, \S1]{SNFBK},
\[
 \HH=\bigoplus_{k\geq0}V^k\mathcal D,
 \qquad \mathcal D=\ker V^*.
\]
On these summands define
\[
 G(V^kx)=(-1)^kV^kx
 \qquad(x\in\mathcal D,\ k\geq0).
\]
Then, $G=G^*=G^{-1}$ and $GVG=-V$. Taking adjoints gives
$GV^*G=-V^*$.

This decomposition identifies $V$ with the unilateral shift on
$\ell^2(\mathbb N_0,\mathcal D)$, where $\mathcal D$ is separable and
infinite-dimensional. Split $\mathcal D$ into four orthogonal
infinite-dimensional subspaces. Each gives a reducing subspace on which
the shift is unitarily equivalent to $V$. Thus,
\[
 V\ueq V\oplus V\oplus V\oplus V.
\]
Together with $V\ueq -V$, this yields
\[
 \lambda V\ueq\diag(\lambda V,\lambda V,-\lambda V,-\lambda V).
\]
Proposition~\ref{prop:fourfold}, applied to $a=\lambda V$, proves
$\lambda V\in\CE$. Taking adjoints in the two unitary equivalences and
applying the same proposition to $a=\lambda V^*$ proves
$\lambda V^*\in\CE$.
\end{proof}

Theorem~\ref{thm:additive} and the following theorem prove
Theorem~\ref{thm:main}.

\begin{theorem}\label{thm:products}
For every $T\in\BH$ and every real $c>\|T\|$, there are $A,B\in\CE$
such that
\[
 T=AB,
 \qquad AA^*=B^*B=cI,
 \qquad \|A\|=\|B\|=\sqrt c.
\]
Consequently,
\[
 \BH=\prod_2\CE=\prod_2\DE.
\]
Every operator is also a product of two operators anticommuting with
involutions. Hence $\ell_\times(\Delta)=2$ for each of the four
classes in Theorem~\ref{thm:main}.
\end{theorem}

\begin{proof}
Apply Theorem~\ref{thm:purefactor} to $C=T/c$, and write $T/c=V^*W$.
Set
\[
 A=\sqrt c\,V^*,
 \qquad B=\sqrt c\,W.
\]
Lemma~\ref{lem:shifts} gives $A,B\in\CE$, while
\[
 AB=T,
 \qquad AA^*=cV^*V=cI,
 \qquad B^*B=cW^*W=cI.
\]
Since $V$ and $W$ are isometries, $\|A\|=\|B\|=\sqrt c$.
Each commutator anticommutes with an involution, is similar to its
negative, and is a difference of idempotents. This gives the other
three product representations.

By Theorem~\ref{thm:additive}, each of the four classes is a proper
subset of $\BH$. Hence one factor does not suffice.
\end{proof}

\medskip

By \cite[Corollary~7.12(a)]{MRZ}, every unitary is a product of two
commutators of idempotents; those factors need not be unitary.
To obtain unitary factors, we use products of bilateral shifts.
Marcoux, Omladi\v{c}, Radjavi and Zhang proved that every unitary is
a product of three bilateral shifts of multiplicity one
\cite[Theorem~5.3]{Tricycles}. The author proved that two suffice
\cite[Theorem~1.1]{Bicycles}. We first represent one such shift as
a commutator of a fixed pair of idempotents.
We use $\mathcal U(\HH)$ for the unitary group.

\begin{proposition}\label{prop:unitary-factors}
There exist an orthogonal projection $P$ and an idempotent $Q$ on
$\HH$ such that $[P,Q]$ is a bilateral shift of
multiplicity one and every
$U\in\mathcal U(\HH)$ has a representation
\[
 U=[R_1PR_1^*,R_1QR_1^*]\,[R_2PR_2^*,R_2QR_2^*]
\]
for some $R_1,R_2\in\mathcal U(\HH)$. In particular,
\[
 \mathcal U(\HH)=\prod_2\bigl(\mathcal U(\HH)\cap\CE\bigr),
\]
and of course one factor does not suffice.
\end{proposition}

\begin{proof}
Identify $\HH$ with $\mathcal L\oplus\mathcal L$, where
$\mathcal L=\ell^2(\mathbb Z)$, and let $We_k=e_{k+1}$ on its
standard orthonormal basis. We need an operator $R$ commuting with
$W$ and satisfying $R^2=W+\frac14I$.

The equation $Wx=-x$ has no nonzero square-summable solution, so
$I+W$ and its adjoint are injective. In the polar decomposition
\[
 I+W=ZF,\qquad F=|I+W|,
\]
$Z$ is unitary and commutes with $F$ and $W$, since $I+W$ is normal.
Moreover,
\[
 Z^2F^2=(I+W)^2=W(2I+W+W^*)=WF^2.
\]
The range of $F^2$ is dense, hence $Z^2=W$.
The norm-convergent binomial series gives
\[
 R=Z\sum_{j=0}^{\infty}\binom{1/2}{j}
          \left(\frac{W^*}{4}\right)^j,
 \qquad RW=WR,\qquad R^2=W+\frac14I.
\]
Set
\[
 P=\begin{pmatrix}I&0\\0&0\end{pmatrix},\qquad
 B=\begin{pmatrix}0&W\\I&0\end{pmatrix}.
\]
Apply \cite[Theorem~2.2]{Drnovsek} with both square roots equal to
$R$. The identities $R^2=W+\tfrac14I$ and $RW=WR$ verify its
hypotheses and give an idempotent $Q$ with $B=[P,Q]$.
We keep this pair $P,Q$ fixed.
The vectors
\[
 f_{2k}=(e_k,0),\qquad f_{2k+1}=(0,e_k)\quad(k\in\mathbb Z)
\]
form an orthonormal basis, and $Bf_j=f_{j+1}$ for every integer $j$.
Thus, $B$ is a bilateral shift of multiplicity one.

By \cite[Theorem~1.1]{Bicycles}, every unitary $U$ can be written as
\[
 U=(R_1BR_1^*)(R_2BR_2^*)
\]
for unitary $R_1,R_2$. Substituting $B=[P,Q]$ proves the stated
representation. Finally, a commutator of idempotents is similar to
its negative, whereas $I$ is not. Thus, one factor does not suffice.
\end{proof}

Thus, for a unitary $U$, the identities $AA^*=B^*B=cI$ hold with
$c=1=\|U\|$. In this case the factors are unitary, so
$A^*A=BB^*=I$ as well.

For differences of idempotents, one can fix the right factor in advance.

\begin{proposition}\label{prop:fixed}
There are a fixed isometry $V\in\CE$ and fixed self-adjoint unitaries
$J,G\in\BH$ such that, for every $T\in\BH$, the operator
\[
 A_T=TV^*-JTV^*J
\]
satisfies
\[
 T=A_TV,
 \qquad JA_TJ=-A_T,
 \qquad GVG=-V.
\]
The map $T\mapsto A_T$ is bounded and linear. Both $A_T$ and $V$
are differences of idempotents.
\end{proposition}

\begin{proof}
Identify $\HH$ with $\bigoplus_{n\geq0}\KK$, where $\KK$ is a
separable infinite-dimensional Hilbert space. Let
\[
 \begin{aligned}
 V(x_0,x_1,x_2,\ldots)&=(0,x_0,x_1,\ldots),\\
 G(x_0,x_1,x_2,\ldots)&=(x_0,-x_1,x_2,\ldots).
 \end{aligned}
\]
Then, $V^*V=I$, $G=G^*=G^{-1}$ and $GVG=-V$.
The isometry $V$ is pure of infinite multiplicity, so
$V\in\CE$ by Lemma~\ref{lem:shifts}.

Put $M=\ran V$ and $N=\ker V^*$. Both spaces are infinite-dimensional
and $\HH=M\oplus N$. Choose a unitary $R:M\to N$ and define
\[
 J(m+n)=Rm+R^*n
 \qquad(m\in M,\ n\in N).
\]
This is a self-adjoint unitary interchanging $M$ and $N$. Therefore,
\[
 V^*JV=0.
\]
For $A_T=TV^*-JTV^*J$, we obtain
\[
 JA_TJ=JTV^*J-TV^*=-A_T
\]
and
\[
 A_TV=TV^*V-JTV^*JV=T.
\]
The formula defines a bounded linear map $T\mapsto A_T$.
Since $A_T$ and $V$ anticommute with involutions, both are differences
of idempotents.
\end{proof}

\section{Nilpotent operators similar to their negatives}\label{sec:nilpotents}

In \cite[Section~2.9]{MRZ}, Marcoux, Radjavi and Zhang ask whether
a nilpotent operator similar to its negative is conjugate to its
negative by an involution. We replace the identity in the upper-right
corner of their block form \cite[Section~2.6]{MRZ} by an operator $C$.
For $D,C\in\mathcal B(\KK)$, put
\begin{equation}\label{eq:weighted-nilpotent}
 M_{D,C}=\begin{pmatrix}0&D&C\\0&0&D\\0&0&0\end{pmatrix}
 \quad\text{on }\HH=\KK\oplus\KK\oplus\KK.
\end{equation}
Throughout this section, $D$ is injective. Then,
$M_{D,C}^3=0\neq M_{D,C}^2$, with $D^2$ in position $(1,3)$ of
$M_{D,C}^2$. We take $D$ and $C$ block diagonal on the same orthogonal
sum. Theorem~\ref{thm:involution-criterion} below determines exactly when an
operator in this family is conjugate to its negative by a bounded
involution. Choosing $C$ compact as well as
$D$ will give compact counterexamples.

\begin{lemma}\label{lem:odd}
Let $\mathcal L$ be a finite-dimensional complex Hilbert space of odd
dimension. For every $A\in\mathcal B(\mathcal L)$ and every unitary
$W\in\mathcal B(\mathcal L)$,
\begin{equation}\label{eq:odd-estimate}
 \|A^2-I\|+\left\|\frac{A+WAW^*}{2}\right\|^2\geq1.
\end{equation}
\end{lemma}

\begin{proof}
Put $\delta=\|A^2-I\|$. If $\delta\geq1$, the inequality is
immediate. Suppose $\delta<1$ and write $A=X+iY$, where $X$ and
$Y$ are self-adjoint. Taking real parts gives
\[
 X^2-Y^2=\frac{A^2+(A^*)^2}{2}\geq(1-\delta)I,
\]
and hence $X^2\geq(1-\delta)I$. This implies that every eigenvalue of $X$
has absolute value at least $a=\sqrt{1-\delta}$.

Since $\dim\mathcal L$ is odd, either the positive or the negative
spectral subspace $\mathcal M$ of $X$ has dimension greater than
$\dim\mathcal L/2$. Choose the corresponding sign
$\eta\in\{1,-1\}$, so that
\[
 \eta\langle Xx,x\rangle\geq a\|x\|^2
 \qquad(x\in\mathcal M).
\]
The spaces $\mathcal M$ and $W\mathcal M$ have a nonzero intersection.
For a unit vector $x$ in that intersection, both $x$ and $W^*x$
belong to $\mathcal M$. Thus,
\[
 \eta\operatorname{Re}\langle(A+WAW^*)x,x\rangle
 =\eta\langle(X+WXW^*)x,x\rangle\geq2a.
\]
It follows that $\|A+WAW^*\|\geq2\sqrt{1-\delta}$, which proves
\eqref{eq:odd-estimate}.
\end{proof}

The unitary blocks in the next theorem are arbitrary. The involution
need not preserve the summands $\KK_n$.

\enlargethispage{\baselineskip}
\begin{theorem}\label{thm:involution-criterion}
Let $\displaystyle\KK=\bigoplus_{n\geq1}\KK_n$, where each $\KK_n$ has finite odd
dimension, and let $U_n$ be unitary on $\KK_n$. Let $(a_n)$ and
$(c_n)$ be bounded sequences of positive numbers. Suppose that, for
some $0<r<1$,
\begin{equation}\label{eq:two-gaps}
 a_{n+1}\leq r a_n,\qquad
 \frac{a_{n+1}^2}{c_{n+1}}\leq r\frac{a_n^2}{c_n}
 \quad\text{for all sufficiently large }n.
\end{equation}
Set
\[
 D=\bigoplus_{n\geq1}a_nU_n,\qquad
 C=\bigoplus_{n\geq1}c_nI_{\KK_n}.
\]
A bounded involution $J$ satisfying $JM_{D,C}J=-M_{D,C}$ exists
if and only if
\begin{equation}\label{eq:bounded-ratio}
 \sup_{n\geq1}\frac{c_n}{a_n}<\infty.
\end{equation}
Moreover, every such involution satisfies, for all sufficiently
large $n$,
\begin{equation}\label{eq:quantitative-involution}
 1\leq\|J\|^2\left(
       \frac{2r(1+r)}{1-r}\frac{a_n}{c_n}
       +\frac{a_n^2}{c_n^2}\right).
\end{equation}

\end{theorem}

\begin{proof}
Suppose first that \eqref{eq:bounded-ratio} holds. The operator
\[
 R=\bigoplus_{n\geq1}\frac{c_n}{a_n}U_n^*
\]
is bounded and satisfies $DR=C$. Direct multiplication gives
\begin{equation}\label{eq:explicit-involution}
 J=\begin{pmatrix}I&0&0\\0&-I&-2R\\0&0&I\end{pmatrix},
 \qquad J^2=I,\qquad JM_{D,C}=-M_{D,C}J.
\end{equation}

Conversely, suppose that $J$ is a bounded involution anticommuting
with $M_{D,C}$, and put $b=\|J\|$. Since $D$ is injective,
\[
 \ker M_{D,C}=\KK\oplus0\oplus0,\qquad
 \ker M_{D,C}^2=\KK\oplus\KK\oplus0.
\]
Both spaces are invariant under $J$. Hence
\[
 J=\begin{pmatrix}
 J_1&J_{12}&J_{13}\\0&J_2&J_{23}\\0&0&J_3
 \end{pmatrix},\qquad J_1^2=J_2^2=J_3^2=I.
\]
Comparing entries gives
\[
 J_1D=-DJ_2,\qquad J_2D=-DJ_3,\qquad
 J_1C+CJ_3=-DJ_{23}-J_{12}D.
\]
In particular,
\begin{equation}\label{eq:weighted-relations}
 J_1D^2=D^2J_3,\qquad
 J_1C+CJ_3=-DJ_{23}-J_{12}D.
\end{equation}
For $C=I$, these are the relations in the proof of
\cite[Proposition~2.8]{MRZ}.

Fix a sufficiently large index $n$, and put $t=a_n$, $s=c_n$ and
$U=U_n$. Use the orthogonal decomposition
\[
 \KK=\KK_-\oplus\KK_n\oplus\KK_+,\qquad
 \KK_-=\bigoplus_{j<n}\KK_j,\quad
 \KK_+=\bigoplus_{j>n}\KK_j.
\]
Write $D=D_-\oplus tU\oplus D_+$ and
$C=C_-\oplus sI\oplus C_+$. Define the bounded operator
\[
 H_+=\bigoplus_{j>n}\frac{a_j^2}{c_j}U_j^2,
 \qquad D_+^2=C_+H_+.
\]
Both sequences $a_j$ and $a_j^2/c_j$ eventually decrease
geometrically by \eqref{eq:two-gaps}. Choose $n$ so large that these
inequalities hold from $n-1$ onward and, for each sequence, the term
at $n-1$ is no greater than any preceding term. Then,
\begin{equation}\label{eq:four-bounds}
 \|D_-^{-1}\|\leq\frac r t,\quad
 \|D_-^{-2}C_-\|\leq\frac{rs}{t^2},\quad
 \|D_+\|\leq rt,\quad
 \|H_+\|\leq\frac{rt^2}{s}.
\end{equation}
The operator $D_-$ is invertible on the finite-dimensional space
$\KK_-$. The second condition in \eqref{eq:two-gaps} makes $H_+$
bounded.

Relative to this decomposition, write
\[
 J_3=\begin{pmatrix} *&X&*\\ Z&A&T\\ *&Y&*\end{pmatrix},
\]
and let $B$ be the middle diagonal entry of $J_1$. Every entry of
$J_1,J_3,J_{12},J_{23}$ has norm at most $b$.

Let $L_-,R_-$ be the entries of $J_{23},J_{12}$ from $\KK_n$ to
$\KK_-$. The corresponding entries in \eqref{eq:weighted-relations},
with the entry of $J_1$ eliminated, give
\[
 sD_-^2X+C_-X(tU)^2
 =-(D_-L_-+R_-tU)(tU)^2.
\]
Multiplication on the left by $s^{-1}D_-^{-2}$ and
\eqref{eq:four-bounds} yield
\[
 (1-r)\|X\|\leq b\frac ts(r+r^2),\qquad
 \|X\|\leq b\frac ts\frac{r(1+r)}{1-r}.
\]

Next consider $T:\KK_+\to\KK_n$. Let $B_+$ be the corresponding
entry of $J_1$, and let $L_+,R_+$ be those of $J_{23},J_{12}$.
Equations \eqref{eq:weighted-relations} give
\[
 B_+D_+^2=(tU)^2T,\qquad
 B_+C_++sT=-tU L_+-R_+D_+.
\]
Multiply the second equality on the right by $H_+$. Since
$C_+H_+=D_+^2$, this gives
\[
 (tU)^2T+sTH_+=-(tU L_++R_+D_+)H_+.
\]
Multiplication on the left by $(tU)^{-2}$ now gives
\[
 (1-r)\|T\|\leq b\frac ts r(1+r),\qquad
 \|T\|\leq b\frac ts\frac{r(1+r)}{1-r}.
\]

The middle diagonal entry of $J_3^2=I$ is
$A^2-I=-ZX-TY$. Thus,
\[
 \|A^2-I\|\leq b(\|X\|+\|T\|)
 \leq b^2\frac ts\frac{2r(1+r)}{1-r}.
\]
Taking middle diagonal entries in \eqref{eq:weighted-relations}
also gives
\[
 B=U^2A(U^*)^2,\qquad \|A+B\|\leq2b\frac ts.
\]
Apply Lemma~\ref{lem:odd} on $\KK_n$ with $W=U^2$. It follows that
\[
 1\leq\|A^2-I\|+\left\|\frac{A+B}{2}\right\|^2
 \leq b^2\left(\frac{2r(1+r)}{1-r}\frac ts+\frac{t^2}{s^2}\right).
\]
This proves \eqref{eq:quantitative-involution}. If $(c_n/a_n)$ were
unbounded, the last expression would be less than one for a suitable
large index. Hence \eqref{eq:bounded-ratio} is necessary.
\end{proof}

The restriction to odd dimensions is needed. On $\mathbb C^2$, put
\[
 U=\begin{pmatrix}1&0\\0&i\end{pmatrix},\qquad
 K=\begin{pmatrix}0&1\\1&0\end{pmatrix}.
\]
Then, $KU^2=-U^2K$, and the self-adjoint unitary
$\diag(-K,-UKU^*,K)$ anticommutes with $M_{aU,cI}$ for every
$a,c>0$. Taking orthogonal sums gives a bounded involution regardless
of the sizes of $c_n/a_n$ in this even-dimensional example.

For classical similarity, an unbounded ratio $c_n/a_n$ can be
compensated by the factor $1+q_n^2$. This is the condition used in
the following construction.

\begin{proposition}\label{prop:weighted-similarity}
Let $(a_n),(c_n)$ be bounded sequences of positive numbers. Let
$d_n\geq3$ be odd integers and choose $q_n\in\mathbb C$ such that
$q_n^{d_n}=1$ and
\begin{equation}\label{eq:weighted-phase}
 \sup_{n\geq1}\frac{c_n}{a_n}|1+q_n^2|<\infty.
\end{equation}
There are unitaries $U_n$ on $\KK_n=\mathbb C^{d_n}$ for which
\[
 D=\bigoplus_{n\geq1}a_nU_n,\qquad
 C=\bigoplus_{n\geq1}c_nI_{\KK_n}
\]
satisfy $M_{D,C}\sim-M_{D,C}$. If the supremum in
\eqref{eq:weighted-phase} is at most $\eta>0$, there is an invertible
$S$ with $S^{-1}M_{D,C}S=-M_{D,C}$ and
\begin{equation}\label{eq:similarity-bound}
 \max\{\|S\|,\|S^{-1}\|\}\leq1+\eta.
\end{equation}
If \eqref{eq:two-gaps} holds and $(c_n/a_n)$ is unbounded, no
bounded involution conjugates $M_{D,C}$ to $-M_{D,C}$. If
$a_n\to0$ and $c_n\to0$, then $M_{D,C}$ is compact.
\end{proposition}

\begin{proof}
On the standard orthonormal basis $e_{n,0},\ldots,e_{n,d_n-1}$, set
\[
 U_ne_{n,j}=q_n^je_{n,j},\qquad
 P_ne_{n,j}=e_{n,j+1\bmod d_n}.
\]
Since $q_n^{d_n}=1$, both operators are unitary and
$U_nP_n=q_nP_nU_n$. Put
\[
 P=\bigoplus_nP_n,\qquad Q=\bigoplus_nq_nI_{\KK_n}.
\]
Then, $P,Q$ are unitary, $DP=QPD$, and $C$ commutes with $D,P,Q$.
Define
\[
 X_n=-\frac{c_n}{a_n}(1+q_n^2)U_n^*P_n,\qquad X=\bigoplus_nX_n.
\]
By \eqref{eq:weighted-phase}, $X$ is bounded and
\[
 DX=-C(I+Q^2)P.
\]
Set
\begin{equation}\label{eq:weighted-S}
 S=\begin{pmatrix}Q^2P&0&0\\0&-QP&X\\0&0&P\end{pmatrix}.
\end{equation}
The three diagonal blocks are unitary and $X$ is bounded, so the
upper triangular operator $S$ is bounded and invertible.
The only possibly nonzero entries of $SM_{D,C}+M_{D,C}S$ are
\[
 \begin{aligned}
 (SM_{D,C}+M_{D,C}S)_{12}&=Q^2PD-DQP=0,\\
 (SM_{D,C}+M_{D,C}S)_{23}&=-QPD+DP=0,\\
 (SM_{D,C}+M_{D,C}S)_{13}&=Q^2PC+DX+CP=0.
 \end{aligned}
\]
Thus, $S^{-1}M_{D,C}S=-M_{D,C}$.

For the norm estimate, let $S_0=\diag(Q^2P,-QP,P)$ and $E=S-S_0$.
The operator $S_0$ is unitary, $\|E\|=\|X\|\leq\eta$, and
$(S_0^*E)^2=0$. Hence
\[
 S^{-1}=(I-S_0^*E)S_0^*,
\]
which proves \eqref{eq:similarity-bound}. Theorem~\ref{thm:involution-criterion}
gives the assertion about involutions. If $a_n,c_n\to0$, the
finite-rank truncations of $D$ and $C$ converge in norm, so every
entry of $M_{D,C}$ is compact.
\end{proof}

We next take $C=|D|^\alpha$. In this case an involution exists
exactly for $\alpha\geq1$.

\begin{corollary}\label{cor:powers}
For each $0<r<1$, there is a compact normal injective operator $D$
such that, for every $\alpha\geq0$,
\[
 M_\alpha=
 \begin{pmatrix}0&D&|D|^\alpha\\0&0&D\\0&0&0\end{pmatrix}
 \quad\text{satisfies}\quad M_\alpha\sim-M_\alpha.
\]
Here $|D|^0=I$. A bounded involution conjugates $M_\alpha$ to
$-M_\alpha$ if and only if $\alpha\geq1$. Moreover, $M_\alpha$ is
compact if and only if $\alpha>0$.
\end{corollary}

\begin{proof}
Choose odd integers $d_n\geq3$ and $d_n$-th roots of unity $q_n$
with $|1+q_n^2|\leq r^n$. Such choices are possible because the roots
of unity of odd order are dense in the unit circle. Let $U_n$ be
as in Proposition~\ref{prop:weighted-similarity}, and set
$D=\bigoplus_nr^nU_n$. Then, $|D|^\alpha=\bigoplus_nr^{\alpha n}I$.
For $a_n=r^n$ and $c_n=r^{\alpha n}$,
\[
 \frac{c_n}{a_n}|1+q_n^2|\leq r^{\alpha n}\leq1.
\]
Proposition~\ref{prop:weighted-similarity} gives $M_\alpha\sim-M_\alpha$
for every $\alpha\geq0$.

If $0\leq\alpha<1$, then $a_n^2/c_n=r^{(2-\alpha)n}$ satisfies
\eqref{eq:two-gaps} with the same $r$, whereas
$c_n/a_n=r^{(\alpha-1)n}$ is unbounded.
Theorem~\ref{thm:involution-criterion} excludes an involution.
If $\alpha\geq1$, the ratio $c_n/a_n$ is bounded and
\eqref{eq:explicit-involution} supplies one. Finally, for
$\alpha>0$ both $D$ and $|D|^\alpha$ are compact. For $\alpha=0$,
the upper-right entry is the identity on an infinite-dimensional
space, so $M_0$ is not compact.
\end{proof}

For $0<\alpha<1$, this gives a compact counterexample. We now impose
summability of every positive power of the singular values.
For $p>0$, let $\mathcal S_p(\HH)$ denote the Schatten $p$-class
of compact operators $T$ on $\HH$ whose singular values
$(s_j(T))_{j\geq1}$ satisfy
\[
 \sum_{j=1}^{\infty}s_j(T)^p<\infty.
\]
For background on singular values and operator ideals, see \cite{Simon}.

\begin{theorem}\label{thm:all-schatten}
For every $\eta>0$, there are
$M\in\displaystyle\bigcap_{p>0}\mathcal S_p(\HH)$ and an invertible $S\in\BH$
such that
\[
 M^3=0\neq M^2,\qquad S^{-1}MS=-M,\qquad
 \max\{\|S\|,\|S^{-1}\|\}\leq1+\eta,
\]
but no bounded involution $J$ satisfies $JMJ=-M$.
\end{theorem}

\begin{proof}
Fix $0<r<1$ and set
\begin{equation}\label{eq:fast-coefficients}
 a_n=r^{n^2+n},\qquad c_n=r^{n^2}\qquad(n\geq1).
\end{equation}
Then,
\[
 \frac{a_{n+1}}{a_n}=r^{2n+2},\qquad
 \frac{a_{n+1}^2/c_{n+1}}{a_n^2/c_n}=r^{2n+3},\qquad
 \frac{c_n}{a_n}=r^{-n}.
\]
Thus, \eqref{eq:two-gaps} holds and $(c_n/a_n)$ is unbounded.

We choose odd $d_n\geq3$ and $d_n$-th roots $q_n$ so that
$$|1+q_n^2|\leq\eta r^n,\qquad d_n\leq L_\eta r^{-n},$$
where $L_\eta$ is independent of $n$.

Now, construct $D,C$ and $M=M_{D,C}$ by
Proposition~\ref{prop:weighted-similarity}. The coefficients tend to
zero, so $M$ is compact. Since
\[
 \frac{c_n}{a_n}|1+q_n^2|\leq\eta,
\]
that proposition gives $S$ with the asserted norm bounds. Theorem~\ref{thm:involution-criterion} excludes every bounded involution, and $M^3=0\neq M^2$ follows from the injectivity of $D$.

Regrouping the orthogonal summands writes $M$ as a direct sum of
matrices $M_n$ on $\KK_n^3$. The range of each $M_n$ is
$\KK_n\oplus\KK_n\oplus\{0\}$, so its rank is $2d_n$. Also,
\[
 \|M_n\|\leq a_n+c_n\leq2c_n.
\]
Consequently, for every $p>0$,
\[
 \sum_j s_j(M)^p
 =\sum_n\sum_j s_j(M_n)^p
 \leq2^{p+1}L_\eta\sum_{n\geq1}r^{pn^2-n}<\infty.
\]
The series converges for every $p>0$. This completes the proof. \end{proof}

\section*{Declarations}

\noindent\textbf{Competing interests.}
The author declares that he has no competing interests.

\smallskip
\noindent\textbf{Funding.}
The author received no funding.

\smallskip
\noindent\textbf{Data availability.}
No datasets were generated during the current study.

\end{document}